\documentclass[11pt,a4paper]{article}

\usepackage{graphicx} % Required for inserting \operatorname{Im}ages
\usepackage[letterpaper,top=2cm,bottom=2cm,left=3cm,right=3cm,marginparwidth=1.75cm]{geometry}
\usepackage{bm}
\usepackage{booktabs}

\usepackage{amsmath}
\usepackage{amssymb,array,longtable} % GKO table

\usepackage[numbers,sort&compress]{natbib}

\usepackage[hidelinks]{hyperref}

\usepackage{caption}
\DeclareCaptionFont{gkotitle}{\fontsize{11}{14.5}\selectfont}
\newcommand{\GKOcaptionnotes}[1]{%
  \endgraf\smallskip
  {\normalfont\fontsize{8}{11}\selectfont\raggedright
   #1\endgraf}%
}
\makeatletter
\renewcommand{\maketitle}{%
  \begingroup
  \centering

  \vspace*{-1.5em}

  \begin{minipage}{1\textwidth}
    \centering

    {\LARGE\bfseries
      \@title
      \par
    }

    \vspace{1.15em}

    {\large
      \@author
      \par
    }

    \vspace{0.85em}

    {\normalsize
      \@date
      \par
    }
  \end{minipage}

  \par
  \vspace{1.4em}
  \endgroup
}
\makeatother

\title{%
Low-Dimensional Reduction of Mean-Field Coupled Dynamical Systems with Units Described by Two Real Variables: Classification and Applications%
}

\author{%
Kai Tokunaga\\[0.55em]
{\small
Department of Complexity Science and Engineering\\
The University of Tokyo, Kashiwa, Chiba 277-8561, Japan\\[0.35em]
\href{mailto:tokunaga.kai25a@c.k.u-tokyo.ac.jp}
{\nolinkurl{tokunaga.kai25a@c.k.u-tokyo.ac.jp}}%
}%
}

\date{20 October 2026}

\begin{document}

\maketitle

\section{abstract}
Low-dimensional reduction is crucial for understanding network systems with many units. In a recent advance, Martens and Dharmavaram developed a theory of exact low-dimensional reduction for a general class of mean-field coupled dynamical systems characterized by a mean-field Lie–Scheffers structure. In this paper, we classify mean-field coupled dynamical systems with a mean-field Lie–Scheffers structure, where each unit is described by two real variables. Our classification is based on González-López, Kamran, and Olver’s local classification of Lie systems on the plane near generic points and yields 28 types of mean-field coupled dynamical systems. For mean-field coupled systems with units described by two real variables, we first present a general theory on classification and low-dimensional reduction based on the mean-field Lie–Scheffers structure, and then discuss applications to populations of coupled oscillators. We show that low-dimensional reduction is possible for a broad class of populations of limit-cycle oscillators, extending the previously reported reducible class of populations of Stuart–Landau oscillators to include many new models.
\section{introduction}
Network systems with many units can be highly complex even when the individual units are simple. Low-dimensional reduction is crucial for analyzing these systems. In particular, many theories of low-dimensional reduction for mean-field coupled dynamical systems have been developed in the study of coupled oscillators. The first major advance was Watanabe–Strogatz (WS) theory \cite{Watanabe1993}. This theory reduces a population of \(N\) identical phase oscillators with sinusoidal coupling to \(N-3\) constants of motion and a three-dimensional system of ODEs. The low-dimensional reduction in WS theory relies on the Riccati structure of sinusoidally coupled phase oscillators \cite{Marvel2009}, and this perspective has led to generalizations of WS theory to populations of units described by complex Riccati equations \cite{Pazo2025}, matrix Riccati equations \cite{Lohe2019}, and Riccati equations on spheres \cite{Lohe2025}. More recently, WS-type reduction theories, which reduce the dynamics to low-dimensional equations and the required number of constants of motion, have also been developed for populations of units governed by quasilinear ODEs \cite{Augustsson2026} and populations of limit-cycle oscillators \cite{tokunaga202608}, which do not have a simple Riccati structure

Recently, Martens and Dharmavaram developed a general theory explaining such WS-type low-dimensional reductions for mean-field coupled dynamical systems with a mean-field Lie–Scheffers structure \cite{martens2026}. Their results show that a mean-field coupled dynamical system with identical units admits exact low-dimensional reduction if the dynamics of each unit form a Lie system \cite{lie1893,Cari_ena_2011} when its mean-field input is replaced by arbitrary time-dependent coefficients. In this paper, we use this result to classify mean-field coupled dynamical systems that admit low-dimensional reduction. According to \cite{martens2026}, a mean-field coupled dynamical system with units described by a Lie system admits low-dimensional reduction. Thus, it is sufficient to consider the classification of Lie systems. Since there are generally infinitely many Lie systems, we restrict the dimension of the phase space of the dynamical systems considered in this paper.

We first consider Lie systems in one real dimension.
These were already classified by Lie himself: their VGL algebras are locally diffeomorphic near generic points to Lie subalgebras of $\operatorname{span}\{\partial_x,x\partial_x,x^2\partial_x\}\simeq\mathfrak{sl}(2,\mathbb{R})$ \cite{Lie1880}.
The most nontrivial and interesting case for mean-field reduction is when each unit is described by a Lie system whose VGL algebra is locally diffeomorphic near generic points to the full algebra $\operatorname{span}\{\partial_x,x\partial_x,x^2\partial_x\}\simeq\mathfrak{sl}(2,\mathbb{R})$.
This is precisely the case in which each unit is described by a Riccati equation, and WS theory corresponds to this case.
What, then, is the classification of Lie systems in two real dimensions?
Gonz\'alez-L\'opez, Kamran, and Olver (GKO) showed that each VGL algebra is locally diffeomorphic near generic points to one of 28 normal forms \cite{GKO1992}.
For mean-field coupled dynamical systems with units described by two real variables, we use this classification to classify those that admit low-dimensional reduction based on the mean-field Lie--Scheffers structure into 28 normal forms.
We then consider applications of this classification to coupled oscillators.

\section{Theory}
\subsection{Problem Formulation}
Let $\bm{x}_i=(x_i,y_i)\in\mathbb{R}^2$, $i=1,\dots,N$, denote the state variables of $N$ units, each described by two real variables.
Defining $\bm{X}=(\bm{x}_1,\dots,\bm{x}_N)$, we can write the mean-field coupled dynamical system with these units as
\begin{equation}
\dot{\bm{X}} = \bm{F}(\bm{X},t),
\end{equation}
\begin{equation}
\bm{F}(\bm{X},t) = \begin{pmatrix}
f(\bm{x}_1,\bm{\mu}(\bm{X}),t) \\
\vdots\\
f(\bm{x}_N,\bm{\mu}(\bm{X}),t)
\end{pmatrix}.
\end{equation}
Here, all units are assumed to obey the same ODE described by a sufficiently smooth function $f:\mathbb{R}^2\times\mathbb{R}^p\times\mathbb{R}\to\mathbb{R}^2$.
The mean-field function $\bm{\mu}:\mathbb{R}^{2N}\to\mathbb{R}^p$ is common to all units.

In \cite{martens2026}, Martens and Dharmavaram defined the mean-field coupled dynamical system (1) to have a mean-field Lie--Scheffers structure if there exists a finite-dimensional real Lie algebra of vector fields $\mathfrak{v}=\operatorname{span}\{L_1,\dots,L_r\}$ such that, for every smooth time-dependent function $\bm{\zeta}:\mathbb{R}\to\mathbb{R}^p$, the vector field $f_x(\bm{x},\zeta(t),t)\partial_x+f_y(\bm{x},\zeta(t),t)\partial_y$ can be written as
\begin{equation}
f_x(\bm{x},\zeta(t),t)\partial_{x}+f_y(\bm{x},\zeta(t),t)\partial_{y} = \sum_{\alpha=1}^{r}a_{\alpha}(t)L_{\alpha}(\bm{x}).
\end{equation}
In \cite{martens2026}, the state of each unit lies in $\mathbb{R}^d$ for general $d$, whereas we restrict our discussion to $d=2$ throughout this paper.

Here, by a finite-dimensional real Lie algebra of vector fields $\mathfrak{v}=\operatorname{span}\{L_1,\dots,L_r\}$, we mean that $L_1,\dots,L_r$ are vector fields on $\mathbb{R}^2$ satisfying
\begin{equation}
[L_{\alpha},L_{\beta}]=\sum_{\gamma=1}^{r}c_{\alpha \beta}^{\gamma}L_{\gamma}
\end{equation}
for all $\alpha,\beta=1,\dots,r$, where $[A,B]=AB-BA$ denotes the Lie bracket.
When the vector field can be expressed as in (3) in terms of a basis $L_1,\dots,L_r$ of a finite-dimensional real Lie algebra, the system
\begin{equation}
\dot{\bm{x}}=f(\bm{x},\zeta(t),t)
\end{equation}
is called a Lie system, and the finite-dimensional real Lie algebra $\mathfrak{v}=\operatorname{span}\{L_1,\dots,L_r\}$ is called a Vessiot--Guldberg Lie (VGL) algebra of the Lie--Scheffers system (5) \cite{lie1893,Cari_ena_2011}.
\subsection{Reduction to fundamental solutions}
In general, the general solution $\bm{x}(t)^{(0)}$ of the Lie--Scheffers system (5) can be expressed in terms of fundamental solutions $\bm{x}^{(1)}(t),\dots,\bm{x}^{(m)}(t)$ using a time-independent superposition function $\Phi:\mathbb{R}^{2m}\times\mathbb{R}^2\to\mathbb{R}^2$:
\begin{equation}
\bm{x}(t)^{(0)}=\Phi(\bm{x}^{(1)}(t),\dots,\bm{x}^{(m)}(t),\bm{\xi}).
\end{equation}
Here, $\bm{\xi}\in\mathbb{R}^2$ is a constant vector determined by the initial value of $\bm{x}(t)^{(0)}$.
Conversely, only Lie--Scheffers systems admit superposition in terms of fundamental solutions as in (6).
Together, these statements constitute the Lie--Scheffers theorem \cite{lie1893,Cari_ena_2011}.

We now briefly describe the reduction to fundamental solutions introduced in \cite{martens2026} for mean-field coupled dynamical systems with a mean-field Lie--Scheffers structure.
Having a mean-field Lie--Scheffers structure means that (5), obtained by replacing the mean-field function for each unit with a general time-dependent function, is a Lie system.
In \cite{martens2026}, the authors show in a subsequent lemma that this replacement is valid when the mean-field function is evaluated along a solution of (1).
Consequently, when (1) has a mean-field Lie--Scheffers structure, applying the Lie--Scheffers theorem to each unit shows that the solution can be expressed using only constants determined by the initial values and $m$ fundamental solutions:
\begin{equation}
\bm{X}(t)^{*}=\big(\Phi(\bm{x}^{(1)}(t),\dots,\bm{x}^{(m)}(t),\bm{\xi}_1),\dots,\Phi(\bm{x}^{(1)}(t),\dots,\bm{x}^{(m)}(t),\bm{\xi}_N)\big).
\end{equation}
The system is then reduced to the following low-dimensional equations for the $m$ fundamental solutions:
\begin{equation}
\begin{pmatrix}
\dot{\bm{x}}^{(1)} \\
\vdots\\
\dot{\bm{x}}^{(m)}
\end{pmatrix} = \begin{pmatrix}
f(\bm{x}^{(1)},\bm{\mu}(\bm{X}^*),t) \\
\vdots\\
f(\bm{x}^{(m)},\bm{\mu}(\bm{X}^*),t)
\end{pmatrix}.
\end{equation}
Here, $\bm{X}^*$ is expressed in terms of $\bm{x}^{(1)}(t),\dots,\bm{x}^{(m)}(t)$ through (7), so (8) is a closed system.

In general, using $m$ fundamental solutions gives a reduction to $2m$ dimensions.
The minimum dimension of the reduced equations obtained by this method is therefore $2m^*$, where $m^*$ is the minimum number of fundamental solutions required.
It was shown in \cite{CARINENA2007} that $m^*$ is the smallest value of $m$ for which, at a generic point $(\bm{x}^{(1)},\dots,\bm{x}^{(m)})$, the only real constants $c_1,\dots,c_r$ satisfying
\begin{equation}
\sum_{\alpha=1}^{r}c_{\alpha}L_{\alpha}(\bm{x}^{(k)})=0,
\quad k=1,\dots,m,
\end{equation}
are $c_1=\dots=c_r=0$.
As we will see below, when (1) has a mean-field Lie--Scheffers structure, it can be reduced to an $r$-dimensional group orbit, where $r$ is the dimension of the VGL algebra.
However, $2m^*$ need not equal $r$ and may be larger.
\subsection{Reduction to equations on a group}
In this paper, we consider reduction to equations on a group in addition to the reduction to fundamental solutions introduced in \cite{martens2026}.
It is known that there exists a matrix Lie algebra isomorphic to the VGL algebra of the Lie system (5) \cite{ado1947,Cari_ena_2011}.
We choose a basis $M_1,\dots,M_r$ of this matrix Lie algebra to satisfy the same commutation relations as $L_1,\dots,L_r$.
Let $G$ be the matrix Lie group generated by this algebra.
Then there exists an action $\Psi:G\times M\to M$ of $G$ on $M$ satisfying the following relation \cite{Cari_ena_2011}:
\begin{equation}
\left.\frac{d}{ds}\right|_{s=0}\Psi(\exp(-sM_{\alpha}),\bm{x})=L_{\alpha}(\bm{x}).
\end{equation}
Using the solution $g$ of the ODE on $G$
\begin{equation}
\dot{g}=A(t)g, \quad g(0)=I, \quad A(t):=-\sum_{\alpha=1}^{r}a_{\alpha}(t)M_{\alpha},
\end{equation}
and a group action $\Psi$ satisfying (10), we can express the general solution of (5) with initial value $\bm{x}_0$ as follows \cite{Cari_ena_2011}:
\begin{equation}
\bm{x}(t)=\Psi(g(t),\bm{x}_0).
\end{equation}

Now consider a mean-field coupled system with a mean-field Lie--Scheffers structure.
As before, we replace the mean-field coupling for each unit with a general time-dependent function and use the general solutions of the individual units to obtain
\begin{equation}
\bm{X}(t)^{*}=\big(\Psi(g(t),\bm{x}_0^{(1)}),\dots,\Psi(g(t),\bm{x}_0^{(N)})\big).
\end{equation}
Here, $\bm{x}_0^{(1)},\dots,\bm{x}_0^{(N)}$ are the initial values of the respective units.
Consequently, the mean-field coupled dynamical system is reduced to the following $r$-dimensional ODE on the group $G$:
\begin{equation}
\dot{g}=A(t)g, \quad g(0)=I, \quad A(t):=-\sum_{\alpha=1}^{r}a_{\alpha}(t,\bm{\mu}(\bm{X}^*))M_{\alpha}.
\end{equation}

We now consider how to obtain an explicit form of the group action $\Psi$ satisfying (10).
For each basis element of the VGL algebra, consider the autonomous system
\begin{equation}
\frac{d}{ds}\psi_{\alpha}^s(\bm{x})=L_{\alpha}(\psi_{\alpha}^s(\bm{x})), \quad \psi_{\alpha}^0(\bm{x})=\bm{x}.
\end{equation}
If $g$ is decomposed as
\begin{equation}
g=\exp(-v_1M_1)\cdots\exp(-v_rM_r),
\end{equation}
we can use the solutions $\psi_{\alpha}^s(\bm{x})$ to construct $\Psi$ as follows \cite{Cari_ena_2011,CARINENA2003}:
\begin{equation}
\Psi(g(t),\bm{x})=(\psi_{1}^{v_1}\circ\cdots\circ\psi_{r}^{v_r})(\bm{x}).
\end{equation}
When the group action is constructed as in (17), substituting it into the ODE for each unit in (1) and requiring the resulting equality to hold identically for general initial values gives an $r$-dimensional system of ODEs for $v_1,\dots,v_r$ equivalent to the ODE (14) on $G$.

Reduction to equations on a group can give a lower-dimensional system than reduction to fundamental solutions.
In some cases, the coordinates $v_1,\dots,v_r$ on the group also have a physical meaning as macroscopic variables, as in WS theory \cite{PIKOVSKY2011}.
In this paper, we use this reduction to equations on a group.
\subsection{Classification of mean-field coupled dynamical systems with a mean-field Lie--Scheffers structure and units described by two real variables}
As we have seen, (1) admits low-dimensional reduction when the system (5), obtained by replacing the mean-field function for each unit with a general time-dependent function, is a Lie system.
For units described by two real variables, we consider the classification of systems for which (5) is a Lie system.
This question has already been answered.
GKO showed that, near generic points, any finite-dimensional Lie algebra of vector fields on the real plane is locally diffeomorphic to one of 28 normal forms \cite{GKO1992}.
However, the VGL algebra associated with a Lie system is not unique, so this does not by itself associate the system with a single one of the 28 normal forms.
To make this correspondence precise, we define the minimal VGL algebra of a Lie system as the VGL algebra of minimum dimension among those associated with the system.
The minimal VGL algebra of the Lie--Scheffers system (5) is then locally diffeomorphic to exactly one of the 28 normal forms.
Even when the phase space is a general real analytic manifold of dimension two, it can be viewed locally as the real plane, so the classification applies to general systems described by two real variables.
However, consistency with global properties, such as the periodicity of the vector fields, must be checked separately.

Table \ref{tab:gko} presents the classification of units in mean-field coupled dynamical systems with a mean-field Lie--Scheffers structure into the 28 GKO normal forms.
A finite-dimensional real Lie algebra $\mathfrak{v}$ is called imprimitive if there exists a family of curves invariant under its action, and primitive otherwise \cite{GKO1992}.
As can be seen from the normal forms, the algebra is imprimitive when the equation for one variable is independent of the other variable.
In the table, $P_1,\dots,P_8$ are primitive types, while $I_1,\dots,I_{11},I_{12}^{(L)},\dots,I_{20}^{(L)}$ are imprimitive types.
The numbering follows that in \cite{BALLESTEROS2015}.
For each normal form, the table also lists the corresponding ODE for the units, the structure of the Lie algebra, the minimum dimension $2m^*$ of the reduction to fundamental solutions, the dimension $r$ of the VGL algebra (the dimension of the equations on the group), and the form of the group action $\Psi$.
The GKO classification does not classify systems by their Lie algebra structures; it gives a classification of dynamical systems into normal forms allowing arbitrary time-dependent coefficients.
As the table shows, some distinct GKO normal forms have the same Lie algebra structure.

\begingroup

% ===== Typesetting settings (unrelated to mathematical coefficients or dimensions) =====
\fontsize{9.5}{13.1}\selectfont % Font size: 9.5 pt; normal baseline spacing: 13.1 pt
\setlength{\tabcolsep}{2pt} % Padding on each side of each cell; the gap between adjacent columns is usually twice this value
\setlength{\arraycolsep}{3.2pt} % Padding on each side of columns in matrices within cells (not the table column widths)
\renewcommand{\arraystretch}{1.15} % Scale factor for the struts determining the minimum row height in tables and arrays
\setlength{\LTleft}{\fill} % Stretchable left margin
\setlength{\LTright}{\fill} % Equal stretchable right margin centers the table
\setlength{\LTcapwidth}{\linewidth} % Caption wrapping width equals the current text width
\setlength{\LTpre}{8pt} % Vertical space before the longtable
\setlength{\LTpost}{8pt} % Vertical space after the longtable

% The [1] in GKOcell specifies one argument. This macro formats cells rather than abbreviating formulas.
% The l column in the array left-aligns each mathematical line within a cell; [t] aligns the top baselines.
% Since raggedright in the column specification does not change the centered alignment within gathered,
% the math macro also uses a single-column, left-aligned array.
\newcommand{\GKOcell}[1]{%
  \(\begin{array}[t]{@{}>{\displaystyle}l@{}}#1\end{array}\)%
}

% Column widths (cell content only, excluding tabcolsep):
% Type 9.7 mm / ODE 43.5 mm / Lie algebra 24.4 mm /
% 2m^* 21.0 mm / r 8.8 mm / Group action 41.4 mm
% Total content width 148.8 mm + intercolumn spacing 20 pt = approximately 155.83 mm.
% The widths have been redistributed to increase intercolumn spacing without exceeding the text width of 155.9 mm.
% @{} removes the outer padding of the table. All six columns, including numerical columns, are left-aligned.
% \\[7pt] at the end of each data row adds space below that row.
% \\*[4pt] at the end of each classification heading adds 4 pt of space and prevents a page break immediately afterward.
% \\[2pt] within a cell adds space between mathematical lines.
\begin{longtable}{@{}
  >{\raggedright\arraybackslash}p{9.7mm}
  >{\raggedright\arraybackslash}p{43.5mm}
  >{\raggedright\arraybackslash}p{24.4mm}
  >{\raggedright\arraybackslash}p{21mm}
  >{\raggedright\arraybackslash}p{8.8mm}
  >{\raggedright\arraybackslash}p{41.4mm}
@{}}
% Within the caption, \endgraf creates a paragraph break (equivalent to \par).
\caption[Classification into normal forms]{Classification into normal forms
\GKOcaptionnotes{%
We use $z=x+iy$ and $\bm u=(x,y)^{\mathsf T}$ as appropriate.
The family parameter $L\ge1$ denotes the number of functions or the degree of a polynomial, and the families depending on $L$ are denoted by $I_{12}^{(L)},\ldots,I_{20}^{(L)}$.
The functions $1,\xi_1,\ldots,\xi_L$ are linearly independent over $\mathbb R$, and $\eta_1,\ldots,\eta_L$ form a fundamental set of solutions of a homogeneous linear ODE of order $L$ with constant coefficients.
Let $D$ denote the differentiation operator $f\mapsto f'$ on $\operatorname{span}_{\mathbb R}\{\eta_1,\ldots,\eta_L\}$, and let $\mathfrak h_2=\mathfrak{aff}(1)$.
The normal forms are understood on the regular domains $y\ne0$ for $P_2$, $x\ne y$ for $I_4$, and $y\ne0$ for $I_5,I_7$.
The coefficients $a_j=a_j(t,\bm\mu)\in\mathbb R$ are common to all units.}}

\label{tab:gko}\label{tab:supp-example}\\
\toprule
Type & ODE for each unit & Lie algebra structure & $2m^*$ & $r$ & Group action $\Psi$\\
\midrule
\endfirsthead
% From the second page onward, repeat only the column headings, without a caption, table number, or continuation notice.
\toprule
Type & ODE for each unit & Lie algebra structure & $2m^*$ & $r$ & Group action $\Psi$\\
\midrule
\endhead
% Use only a horizontal rule at the bottom of each page, without a caption.
\bottomrule
\endfoot
\bottomrule
\endlastfoot
% Place common conditions on coefficients and group parameters inside the table rather than in the caption.
\multicolumn{6}{@{}l}{\textit{Primitive}: $P_1,\ldots,P_{8}$}\\*[4pt]
% P1
$P_{1}$
& \GKOcell{\dot z=a_1+ia_2+(\alpha-i)a_3z,\\\alpha\ge0}
& \GKOcell{A_\alpha\simeq\\\mathbb R\ltimes\mathbb R^2}
& \GKOcell{4}
& \GKOcell{3}
& \GKOcell{z\mapsto e^{(\alpha-i)s}z+\beta,\\s\in\mathbb R,\ \beta\in\mathbb C} \\[7pt]
% P2
$P_{2}$
& \GKOcell{\dot z=a_1+a_2z+a_3z^2}
& \GKOcell{\mathfrak{sl}(2,\mathbb R)}
& \GKOcell{4}
& \GKOcell{3}
& \GKOcell{z\mapsto\frac{az+b}{cz+d},\\a,b,c,d\in\mathbb R,\\ad-bc=1} \\[7pt]
% P3
$P_{3}$
& \GKOcell{\begin{aligned}[t]\dot z&=a_2+ia_3-ia_1z\\&\quad+(a_2-ia_3)z^2\end{aligned}}
& \GKOcell{\mathfrak{so}(3)}
& \GKOcell{4}
& \GKOcell{3}
& \GKOcell{z\mapsto\frac{az+b}{-\bar b z+\bar a},\\a,b\in\mathbb C,\\|a|^2+|b|^2=1} \\[7pt]
% P4
$P_{4}$
& \GKOcell{\dot z=a_1+ia_2+(a_3-ia_4)z}
& \GKOcell{\mathbb R^2\ltimes\mathbb R^2}
& \GKOcell{4}
& \GKOcell{4}
& \GKOcell{z\mapsto az+b,\\a\in\mathbb C\setminus\{0\},\\b\in\mathbb C} \\[7pt]
% P5
$P_{5}$
& \GKOcell{\dot{\bm u}=M\bm u+\bm b,\\[2pt]
M=\begin{pmatrix}a_3&a_4\\a_5&-a_3\end{pmatrix},\\[2pt]
\bm b=(a_1,a_2)^{\mathsf T}}
& \GKOcell{\mathfrak{sl}(2,\mathbb R)\\{}\ltimes\mathbb R^2}
& \GKOcell{6}
& \GKOcell{5}
& \GKOcell{\bm u\mapsto A\bm u+\bm\beta,\\[2pt]
A=\begin{pmatrix}a&b\\c&d\end{pmatrix}\in\mathbb R^{2\times2},\\[2pt]
\bm\beta\in\mathbb R^2,\quad ad-bc=1} \\[7pt]
% P6
$P_{6}$
& \GKOcell{\dot{\bm u}=M\bm u+\bm b,\\[2pt]
M=\begin{pmatrix}a_3&a_4\\a_5&a_6\end{pmatrix},\\[2pt]
\bm b=(a_1,a_2)^{\mathsf T}}
& \GKOcell{\mathfrak{gl}(2,\mathbb R)\\{}\ltimes\mathbb R^2}
& \GKOcell{6}
& \GKOcell{6}
& \GKOcell{\bm u\mapsto A\bm u+\bm\beta,\\[2pt]
A=\begin{pmatrix}a&b\\c&d\end{pmatrix}\in\mathbb R^{2\times2},\\[2pt]
\bm\beta\in\mathbb R^2,\quad ad-bc>0} \\[7pt]
% P7
$P_{7}$
& \GKOcell{\begin{aligned}[t]\dot z&=a_1+ia_2+(a_3-ia_4)z\\&\quad+(a_5-ia_6)z^2\end{aligned}}
& \GKOcell{\mathfrak{so}(3,1)\\\simeq\mathfrak{sl}(2,\mathbb C)_{\mathbb R}}
& \GKOcell{6}
& \GKOcell{6}
& \GKOcell{z\mapsto\frac{az+b}{cz+d},\\a,b,c,d\in\mathbb C,\\ad-bc=1} \\[7pt]
% P8
$P_{8}$
& \GKOcell{\dot{\bm u}=\bm b+M\bm u+(\bm c^{\mathsf T}\bm u)\bm u,\\[2pt]M=\begin{pmatrix}a_3&a_4\\a_5&a_6\end{pmatrix},\\[2pt]
\bm b=(a_1,a_2)^{\mathsf T},\\
\bm c=(a_7,a_8)^{\mathsf T}}
& \GKOcell{\mathfrak{sl}(3,\mathbb R)}
& \GKOcell{8}
& \GKOcell{8}
& \GKOcell{\bm u\mapsto\frac{A\bm u+\bm\beta}{\bm\gamma^{\mathsf T}\bm u+d},\\[2pt]
\det\begin{pmatrix}A&\bm\beta\\\bm\gamma^{\mathsf T}&d\end{pmatrix}=1,\\[2pt]
A\in\mathbb R^{2\times2},\quad d\in\mathbb R,\\
\bm\beta,\bm\gamma\in\mathbb R^2} \\[7pt]
\addlinespace[4pt]
\multicolumn{6}{@{}l}{\textit{Imprimitive}: $I_1,\ldots,I_{11},\ I_{12}^{(L)},\ldots,I_{20}^{(L)}$}\\*[4pt]
% I1
$I_{1}$
& \GKOcell{\dot x=a_1,\\\dot y=0}
& \GKOcell{\mathbb R}
& \GKOcell{2}
& \GKOcell{1}
& \GKOcell{(x,y)\mapsto(x+u,y)} \\[7pt]
% I2
$I_{2}$
& \GKOcell{\dot x=a_1+a_2x,\\\dot y=0}
& \GKOcell{\mathfrak h_2}
& \GKOcell{4}
& \GKOcell{2}
& \GKOcell{(x,y)\mapsto(\lambda x+u,y),\\\lambda>0} \\[7pt]
% I3
$I_{3}$
& \GKOcell{\dot x=a_1+a_2x+a_3x^2,\\\dot y=0}
& \GKOcell{\mathfrak{sl}(2,\mathbb R)\\\text{type I}}
& \GKOcell{6}
& \GKOcell{3}
& \GKOcell{(x,y)\mapsto\left(\frac{ax+b}{cx+d},y\right),\\ad-bc=1} \\[7pt]
% I4
$I_{4}$
& \GKOcell{\begin{aligned}[t]\dot x&=a_1+a_2x+a_3x^2,\\\dot y&=a_1+a_2y+a_3y^2\end{aligned}}
& \GKOcell{\mathfrak{sl}(2,\mathbb R)\\\text{type II}}
& \GKOcell{4}
& \GKOcell{3}
& \GKOcell{(x,y)\mapsto\\
\left(\frac{ax+b}{cx+d},\frac{ay+b}{cy+d}\right),\\ad-bc=1} \\[7pt]
% I5
$I_{5}$
& \GKOcell{\begin{aligned}[t]\dot x&=a_1+2a_2x+a_3x^2,\\\dot y&=a_2y+a_3xy\end{aligned}}
& \GKOcell{\mathfrak{sl}(2,\mathbb R)\\\text{type III}}
& \GKOcell{4}
& \GKOcell{3}
& \GKOcell{(x,y)\mapsto\\
\left(\frac{ax+b}{cx+d},\frac{y}{cx+d}\right),\\ad-bc=1} \\[7pt]
% I6
$I_{6}$
& \GKOcell{\dot x=a_1+a_3x+a_4x^2,\\\dot y=a_2}
& \GKOcell{\mathfrak{gl}(2,\mathbb R)\\\text{type I}}
& \GKOcell{6}
& \GKOcell{4}
& \GKOcell{(x,y)\mapsto\\
\left(\frac{ax+b}{cx+d},y+u\right),\\ad-bc=1} \\[7pt]
% I7
$I_{7}$
& \GKOcell{\begin{aligned}[t]\dot x&=a_1+a_3x+a_4x^2,\\\dot y&=a_2y+a_4xy\end{aligned}}
& \GKOcell{\mathfrak{gl}(2,\mathbb R)\\\text{type II}}
& \GKOcell{4}
& \GKOcell{4}
& \GKOcell{(x,y)\mapsto\\
\left(\frac{ax+b}{cx+d},\frac{y}{cx+d}\right),\\ad-bc>0} \\[7pt]
% I8
$I_{8}$
& \GKOcell{\begin{aligned}[t]\dot x&=a_1+a_3x,\\\dot y&=a_2+\alpha a_3y\end{aligned}\\0<|\alpha|\le1}
& \GKOcell{B_\alpha\simeq\\\mathbb R\ltimes\mathbb R^2}
& \GKOcell{4}
& \GKOcell{3}
& \GKOcell{(x,y)\mapsto\\
(\lambda x+u,\lambda^\alpha y+v),\\\lambda>0} \\[7pt]
% I9
$I_{9}$
& \GKOcell{\begin{aligned}[t]\dot x&=a_1+a_3x,\\\dot y&=a_2+a_4y\end{aligned}}
& \GKOcell{\mathfrak h_2\oplus\mathfrak h_2}
& \GKOcell{4}
& \GKOcell{4}
& \GKOcell{(x,y)\mapsto\\
(\lambda x+u,\nu y+v),\\\lambda,\nu>0} \\[7pt]
% I10
$I_{10}$
& \GKOcell{\begin{aligned}[t]\dot x&=a_1+a_3x+a_5x^2,\\\dot y&=a_2+a_4y\end{aligned}}
& \GKOcell{\mathfrak{sl}(2,\mathbb R)\\{}\oplus\mathfrak h_2}
& \GKOcell{6}
& \GKOcell{5}
& \GKOcell{(x,y)\mapsto\\
\left(\frac{ax+b}{cx+d},\lambda y+u\right),\\ad-bc=1,\quad\lambda>0} \\[7pt]
% I11
$I_{11}$
& \GKOcell{\begin{aligned}[t]\dot x&=a_1+a_3x+a_5x^2,\\\dot y&=a_2+a_4y+a_6y^2\end{aligned}}
& \GKOcell{\mathfrak{sl}(2,\mathbb R)\\{}\oplus\mathfrak{sl}(2,\mathbb R)}
& \GKOcell{6}
& \GKOcell{6}
& \GKOcell{(x,y)\mapsto\\
\left(\frac{ax+b}{cx+d},\frac{ey+f}{hy+\ell}\right),\\ad-bc=1,\quad e\ell-fh=1} \\[7pt]
% I12
$I_{12}^{(L)}$
& \GKOcell{\begin{aligned}[t]\dot x&=0,\\\dot y&=a_1+\sum_{j=1}^L a_{j+1}\xi_j(x)\end{aligned}}
& \GKOcell{\mathbb R^{L+1}}
& \GKOcell{2L+2}
& \GKOcell{L+1}
& \GKOcell{(x,y)\mapsto\\
\bigl(x,y+h_0+\sum_{j=1}^L h_j\xi_j(x)\bigr)} \\[7pt]
% I13
$I_{13}^{(L)}$
& \GKOcell{\begin{aligned}[t]\dot x&=0,\\\dot y&=a_1+a_2y+\sum_{j=1}^L a_{j+2}\xi_j(x)\end{aligned}}
& \GKOcell{\mathbb R\ltimes\mathbb R^{L+1}}
& \GKOcell{2L+4}
& \GKOcell{L+2}
& \GKOcell{(x,y)\mapsto\\
\bigl(x,\lambda y+h_0+\sum_{j=1}^L h_j\xi_j(x)\bigr)\\\lambda>0} \\[7pt]
% I14
$I_{14}^{(L)}$
& \GKOcell{\dot x=a_1,\\\dot y=\sum_{j=1}^L a_{j+1}\eta_j(x)}
& \GKOcell{\mathbb R\ltimes_D\mathbb R^L}
& \GKOcell{2L}
& \GKOcell{L+1}
& \GKOcell{(x,y)\mapsto\\
\bigl(x+u,y+\sum_{j=1}^L h_j\eta_j(x)\bigr)} \\[7pt]
% I15
$I_{15}^{(L)}$
& \GKOcell{\begin{aligned}[t]\dot x&=a_1,\\\dot y&=a_2y+\sum_{j=1}^L a_{j+2}\eta_j(x)\end{aligned}}
& \GKOcell{\mathbb R^2\\{}\ltimes_{(D,-\mathrm{id})}\mathbb R^L}
& \GKOcell{2L+2}
& \GKOcell{L+2}
& \GKOcell{(x,y)\mapsto\\
\bigl(x+u,\lambda y+\sum_{j=1}^L h_j\eta_j(x)\bigr)\\\lambda>0} \\[7pt]
% I16
$I_{16}^{(L)}$
& \GKOcell{\begin{aligned}[t]\dot x&=a_1+a_2x,\\\dot y&=\alpha a_2y+\sum_{j=0}^L a_{j+3}x^j\end{aligned}\\\alpha\in\mathbb R}
& \GKOcell{C_\alpha^L\simeq\\\mathfrak h_2\ltimes\mathbb R^{L+1}}
& \GKOcell{2L+2}
& \GKOcell{L+3}
& \GKOcell{(x,y)\mapsto\\
\bigl(\lambda x+u,\lambda^\alpha y+\sum_{j=0}^L h_jx^j\bigr)\\\lambda>0} \\[7pt]
% I17
$I_{17}^{(L)}$
& \GKOcell{\begin{aligned}[t]\dot x&=a_1+a_2x,\\\dot y&=a_2(Ly+x^L)\\&\quad+\sum_{j=0}^{L-1}a_{j+3}x^j\end{aligned}}
& \GKOcell{\mathbb R\ltimes\\(\mathbb R\ltimes\mathbb R^L)}
& \GKOcell{2\max(2,L)}
& \GKOcell{L+2}
& \GKOcell{(x,y)\mapsto\\
\begin{aligned}[t]
\bigl(&\lambda x+u,\\
&\lambda^Ly+\lambda^L\log\lambda\,x^L\\
&\quad+\sum_{j=0}^{L-1}h_jx^j\bigr)
\end{aligned}\\\lambda>0} \\[7pt]
% I18
$I_{18}^{(L)}$
& \GKOcell{\begin{aligned}[t]\dot x&=a_1+a_2x,\\\dot y&=a_3y+\sum_{j=0}^L a_{j+4}x^j\end{aligned}}
& \GKOcell{(\mathfrak h_2\oplus\mathbb R)\\{}\ltimes\mathbb R^{L+1}}
& \GKOcell{2L+4}
& \GKOcell{L+4}
& \GKOcell{(x,y)\mapsto\\
\bigl(\lambda x+u,\nu y+\sum_{j=0}^L h_jx^j\bigr)\\\lambda,\nu>0} \\[7pt]
% I19
$I_{19}^{(L)}$
& \GKOcell{\begin{aligned}[t]\dot x&=a_1+2a_2x+a_3x^2,\\\dot y&=La_2y+La_3xy\\&\quad+\sum_{j=0}^L a_{j+4}x^j\end{aligned}}
& \GKOcell{\mathfrak{sl}(2,\mathbb R)\\{}\ltimes\mathbb R^{L+1}}
& \GKOcell{2\max\!(3,L\!+\!1)}
& \GKOcell{L+4}
& \GKOcell{(x,y)\mapsto\\
\left(\frac{ax+b}{cx+d},\frac{y+\sum\limits_{j=0}^L h_jx^j}{(cx+d)^L}\right),\\ad-bc=1} \\[7pt]
% I20
$I_{20}^{(L)}$
& \GKOcell{\begin{aligned}[t]\dot x&=a_1+a_2x+a_4x^2,\\\dot y&=a_3y+La_4xy\\&\quad+\sum_{j=0}^L a_{j+5}x^j\end{aligned}}
& \GKOcell{\mathfrak{gl}(2,\mathbb R)\\{}\ltimes\mathbb R^{L+1}}
& \GKOcell{2L+4}
& \GKOcell{L+5}
& \GKOcell{(x,y)\mapsto\\
\left(\frac{ax+b}{cx+d},\frac{\lambda y+\sum\limits_{j=0}^L h_jx^j}{(cx+d)^L}\right),\\ad-bc=1,\quad\lambda>0} \\[7pt]
\end{longtable}
\endgroup

\section{Applications to coupled oscillators}
We now consider applications to coupled oscillators.
In the study of coupled oscillators, low-dimensional reduction for populations of phase oscillators is well known, but relatively few results are available for more general populations of limit-cycle oscillators that include amplitude dynamics.
Below, we show that a broad class of populations of limit-cycle oscillators admits low-dimensional reduction.
For simplicity, we restrict our attention to systems in which the phase obeys a closed equation of Riccati type.
These classes are all imprimitive and are obtained from $I_3,I_4,I_5,I_6,I_7,I_{10},I_{11},I_{19}^{(L)},I_{20}^{(L)}$ by setting $x=\tan\frac{k\theta}{2}$.
For the amplitude coordinate, we set $y=\rho$ for $I_3,I_4,I_6,I_{10},I_{11}$, $y=\rho(1+x^2)^{1/2}$ for $I_5,I_7$, and $y=\rho(1+x^2)^{L/2}$ for $I_{19}^{(L)},I_{20}^{(L)}$.
For $I_{19}^{(L)}$ and $I_{20}^{(L)}$, at least one of $k$ and $L$ must be even for the phase variable to be well defined globally.
If we further restrict attention to systems that can have a stable limit cycle, the possible types are $I_7,I_{10},I_{11},I_{20}^{(L)}$.
Table \ref{tab:phase-amplitude} summarizes the corresponding ODEs for the units in phase--amplitude coordinates.

Among these types, we focus on $I_{20}^{(L)}$ because it allows phase-dependent amplitude dynamics and has considerable generality.
For $I_{20}^{(L)}$, we set $\omega(\bm{\mu}(\bm{X}),t)=\omega$ and $\lambda(\bm{\mu}(\bm{X}),t)=\lambda$, and assume that all other terms represent interactions only.
In the absence of interactions, the system then becomes
\begin{equation}
\begin{split}
\dot{\theta}_i &= \omega,\\
\dot{\rho}_i &= -\lambda\rho_i.
\end{split}
\end{equation}
For $\lambda>0$, this can be interpreted as an unperturbed limit-cycle oscillator expressed in phase--amplitude coordinates.
When interactions are included, the resulting system can be interpreted as a phase--amplitude reduction through first order, $O(\rho)$, in the deviation from the limit cycle \cite{Wilson2016,Shirasaka2020}.
However, when the system is interpreted as a phase--amplitude oscillator in this way, the forms of its response functions must be regarded as strongly constrained.

A more natural interpretation is to regard $r_i$ as the physical amplitude, with $\rho_i=1/r_i^p$ for a positive integer $p$.
The equations in phase--amplitude coordinates then become
\begin{equation}
\begin{split}
\dot{\theta}_i
&=\omega(\bm{\mu}(\bm{X}),t)
+\operatorname{Im}\left(
H(\bm{\mu}(\bm{X}),t)e^{ik\theta_i}
\right),\\
\dot{r}_i
&=\left[
\lambda(\bm{\mu}(\bm{X}),t)
-\frac{Lk}{2}\operatorname{Re}\left(
H(\bm{\mu}(\bm{X}),t)e^{ik\theta_i}
\right)
\right]\frac{r_i}{p}\\
&\quad-\frac{1}{p}\sum_{j=0}^{L/2}
\operatorname{Re}\left(
G_j^{(L)}(\bm{\mu}(\bm{X}),t)
e^{\frac{i(L-2j)k\theta_i}{2}}
\right)r_i^{p+1}.
\end{split}
\end{equation}
Here, we consider only even values of $L$.
Further setting $z_i=r_ie^{i\theta_i}$ gives
\begin{equation}
\begin{split}
\dot{z}_i
&=\left[
\frac{\lambda(\bm{\mu}(\bm{X}),t)}{p}
+i\omega(\bm{\mu}(\bm{X}),t)
\right]z_i\\
&\quad+\left[
\frac{2p-Lk}{4p}
H(\bm{\mu}(\bm{X}),t)
\left(\frac{z_i}{|z_i|}\right)^k
\right.\\
&\qquad\left.
-\frac{Lk+2p}{4p}
\overline{H}(\bm{\mu}(\bm{X}),t)
\left(\frac{\overline{z}_i}{|z_i|}\right)^k
\right]z_i\\
&\quad-\frac{1}{p}\sum_{j=0}^{L/2}
\operatorname{Re}\left(
G_j^{(L)}(\bm{\mu}(\bm{X}),t)
\left(\frac{z_i}{|z_i|}\right)^{\frac{(L-2j)k}{2}}
\right)z_i|z_i|^p.
\end{split}
\end{equation}
This is precisely a Poincar\'e oscillator \cite{glass1989} for $p=1$ and a Stuart--Landau oscillator \cite{landau1944,Stuart_1960} for $p=2$.
For $p\ge3$, it is a generalization of the Stuart--Landau oscillator with higher-order saturation.

In many cases, the vector field loses smoothness at the origin.
These models are therefore often unsuitable for describing bifurcations in which oscillations emerge from a state of zero amplitude.
However, there are also cases in which the vector field remains smooth at the origin.
We consider such cases for $p=2$.
Since it suffices for the vector field to be polynomial in $z_i$ and $\overline{z}_i$, all cases that are smooth at the origin for $p=2$ are given by the following form, with $p=2,k=2,L=2$:
\begin{equation}
\begin{split}
\dot{z}_i
&=\left[
\frac{\lambda(\bm{\mu}(\bm{X}),t)}{2}
+i\omega(\bm{\mu}(\bm{X}),t)
\right]z_i\\
&\quad-\frac{1}{2}\operatorname{Re}\left(
G_1^{(2)}(\bm{\mu}(\bm{X}),t)
\right)z_i|z_i|^2\\
&\quad-\overline{H}(\bm{\mu}(\bm{X}),t)\overline{z_i}
-\frac{1}{2}\operatorname{Re}\left(
G_0^{(2)}(\bm{\mu}(\bm{X}),t)z_i^2
\right)z_i.
\end{split}
\end{equation}
This is precisely the population of Stuart--Landau oscillators with coupling beyond the coefficients that was shown to admit low-dimensional reduction in \cite{tokunaga202608} (the case of coupling through the coefficients belongs to $I_{15}^{(L)}$).
Because this model has second-harmonic coupling, it typically exhibits clustering.
First-harmonic coupling can also be considered, although smoothness is lost at the origin.
For example, for populations of Stuart--Landau oscillators, we can consider the following model with first-harmonic coupling $(p=2,k=1,L=4)$:
\begin{equation}
\begin{split}
\dot{z}_i
&=\left[
\frac{\lambda(\bm{\mu}(\bm{X}),t)}{2}
+i\omega(\bm{\mu}(\bm{X}),t)
\right]z_i\\
&\quad-\frac{1}{2}\operatorname{Re}\left(
G_2^{(4)}(\bm{\mu}(\bm{X}),t)
\right)z_i|z_i|^2
-\overline{H}(\bm{\mu}(\bm{X}),t)|z_i|.
\end{split}
\end{equation}
For populations of Poincar\'e oscillators, we obtain the following model $(p=1,k=1,L=2)$:
\begin{equation}
\begin{split}
\dot{z}_i
&=\left[
\lambda(\bm{\mu}(\bm{X}),t)
+i\omega(\bm{\mu}(\bm{X}),t)
\right]z_i\\
&\quad-\operatorname{Re}\left(
G_1^{(2)}(\bm{\mu}(\bm{X}),t)
\right)z_i|z_i|
-\overline{H}(\bm{\mu}(\bm{X}),t)|z_i|.
\end{split}
\end{equation}
Equations (19) and (20) also describe a variety of other models for populations of limit-cycle oscillators.

\begingroup
\fontsize{9.5}{13.1}\selectfont
\setlength{\tabcolsep}{2pt}
\setlength{\arraycolsep}{3.2pt}
\renewcommand{\arraystretch}{1.15}
\setlength{\LTleft}{\fill}
\setlength{\LTright}{\fill}
% Set the width of Table 2 to suit the equations rather than extending it to the full text width.
% Changing this 100 mm value adjusts both the table width and the caption width.
\newcommand{\PATableWidth}{100mm}
\setlength{\LTcapwidth}{\PATableWidth}
\setlength{\LTpre}{8pt}
\setlength{\LTpost}{8pt}
\newcommand{\PAcell}[1]{%
  \(\begin{array}[t]{@{}>{\displaystyle}l@{}}#1\end{array}\)%
}
% The total table width is 100 mm. The Type column is 12 mm wide, with the remaining width allocated to the ODE column.
% Equal LTleft and LTright values center the table within the text area.
% Cell contents remain left-aligned. The caption is also aligned with the left edge of the table.
\begin{longtable}{@{}
  >{\raggedright\arraybackslash}p{12mm}
  >{\raggedright\arraybackslash}p{\dimexpr\PATableWidth-12mm-2\tabcolsep\relax}
@{}}
\caption[Populations of limit-cycle oscillators]{Populations of limit-cycle oscillators
}
\label{tab:phase-amplitude}\\
\toprule
Type & ODE for each unit in phase--amplitude coordinates $(\theta,\rho)$\\
\midrule
\endfirsthead
\toprule
Type & ODE for each unit in phase--amplitude coordinates $(\theta,\rho)$\\
\midrule
\endhead
\bottomrule
\endfoot
\bottomrule
\endlastfoot
% Since this is a short table with four rows, prevent page breaks between the data rows.
$I_7$
& \PAcell{\begin{aligned}[t]
\dot\theta&=\omega+\operatorname{Im}(He^{ik\theta}),\\[3pt]
\dot\rho&=\left[-\lambda+\frac{k}{2}\operatorname{Re}(He^{ik\theta})\right]\rho
\end{aligned}}\\*[7pt]
$I_{10}$
& \PAcell{\begin{aligned}[t]
\dot\theta&=\omega+\operatorname{Im}(He^{ik\theta}),\\[3pt]
\dot\rho&=-\lambda\rho+\gamma_0
\end{aligned}}\\*[7pt]
$I_{11}$
& \PAcell{\begin{aligned}[t]
\dot\theta&=\omega+\operatorname{Im}(He^{ik\theta}),\\[3pt]
\dot\rho&=\gamma_0-\lambda\rho+\gamma_2\rho^2
\end{aligned}}\\*[7pt]
$I_{20}^{(L)}$
& \PAcell{\begin{aligned}[t]
\dot\theta&=\omega+\operatorname{Im}(He^{ik\theta}),\\[3pt]
\dot\rho&=\left[-\lambda+\frac{Lk}{2}\operatorname{Re}(He^{ik\theta})\right]\rho
{}+\sum_{j=0}^{\lfloor L/2\rfloor}
\operatorname{Re}\!\left(G_j^{(L)}e^{\frac{i(L-2j)k\theta}{2}}\right)
\end{aligned}}\\[7pt]
\end{longtable}
\endgroup
\section{Conclusion}
In this paper, we considered mean-field coupled dynamical systems with identical units described by two real variables and classified those that admit low-dimensional reduction based on the mean-field Lie--Scheffers structure into 28 types using the GKO classification of Lie systems.
We first showed that, in addition to reduction based on the superposition of fundamental solutions, the mean-field Lie--Scheffers structure allows reduction to equations on a Lie group using group actions.
We then classified the systems using the GKO classification and finally applied the results to the low-dimensional reduction of populations of limit-cycle oscillators.
The populations of limit-cycle oscillators identified in this study admit low-dimensional reduction and substantially generalize the reducible populations of Stuart--Landau oscillators identified in previous work \cite{tokunaga202608}.
These models include not only the second-harmonic coupling considered in \cite{tokunaga202608} but also first-harmonic coupling.

An important direction for future work related to the results of this paper is to investigate mean-field coupled dynamical systems with units of higher dimension.
For example, the case of two complex variables can be addressed immediately, since a classification into 23 normal forms is already known \cite{Lie1880}.
For three or more real variables, however, a complete classification is not yet known, and the number of types is expected to be very large.
For example, even among systems in three real variables with three-dimensional VGL algebras, 31 types are known \cite{IBRAGIMOV2016}.
It therefore seems reasonable to focus on interesting classes in one or two variables and generalize them to higher dimensions.
For example, low-dimensional reduction has already been obtained for mean-field coupled systems with units described by Riccati systems generalized to higher-dimensional spheres \cite{Lohe2025}.
Recently, Stuart--Landau oscillators have also been generalized to higher dimensions \cite{Gogoi2024}.
It would therefore be interesting to study low-dimensional reduction for higher-dimensional generalizations of the populations of limit-cycle oscillators of type $I_{20}^{(L)}$ identified in this study.

\section*{Acknowledgments}

The author thanks Hiroshi Kori for valuable discussions. This study was supported by the WINGS-FMSP program at the University of Tokyo.

% ============================================================
% References
% ============================================================
\bibliographystyle{unsrtnat}
\bibliography{Ref}

\end{document}